# Computational Models based on Synchronized Oscillators for Solving Combinatorial Optimization Problems

Antik Mallick[1#], Mohammad Khairul Bashar[1#], Zongli Lin[1], Nikhil Shukla[1*]

[1]Department of Electrical & Computer Engineering, University of Virginia, Charlottesville, VA- 22904, USA

[#]Equal contribution

*e-mail: ns6pf@virginia.edu

**Abstract:**

The equivalence between the natural minimization of energy in a dynamical system and the minimization of an objective function characterizing a combinatorial optimization problem offers a promising approach to designing dynamical system inspired computational models and solvers for such problems. For instance, the ground state energy of coupled electronic oscillators, under second harmonic injection, can be directly mapped to the optimal solution of the Maximum Cut problem. However, prior work has focused on a limited set of such problems. Therefore, in this work, we formulate computing models based on synchronized oscillator dynamics for a broad spectrum of combinatorial optimization problems ranging from the Max-K-Cut (the general version of the Maximum Cut problem) to the Traveling Salesman Problem. We show that synchronized oscillator dynamics can be engineered to solve these different combinatorial optimization problems by appropriately designing the coupling function and the external injection to the oscillators. Our work marks a step forward towards expanding the functionalities of oscillator-based analog accelerators, and furthers the scope of dynamical system solvers for combinatorial optimization problems.

**Introduction:**

The rich spectrum of spatiotemporal properties offered by dynamical systems such as synchronized oscillators can facilitate novel physics-inspired pathways to solving computationally challenging problems [1-9]. As a case in point, the phase dynamics of coupled oscillators, under second harmonic injection, can be directly mapped to the solution of the Ising Hamiltonian $H = -\sum_{I,j,i<j}^{N} J_{ij}s_is_j$ (the Zeeman term has been neglected here) allowing for the creation of oscillator-based Ising machines [10] (described in the following section). Subsequently, such oscillator Ising machines can be used to solve many computationally intractable problems that can be expressed efficiently in terms of minimizing the Ising Hamiltonian. For instance, many prior theoretical and experimental demonstrations of Ising machines in literature have been shown to solve the archetypal Maximum Cut (MaxCut) problem [11-20] – the challenge of dividing the nodes of a graph into two sets such that the number of common edges among the two sets is as large as possible.

However, prior work [21-24] has primarily focused on formulating oscillator-based computational models for a relatively small subset of combinatorial problems. For instance, while oscillator-based models are available for the MaxCut problem – a special instance of the Max-K-Cut (with K=2), no such models have been developed, to the best of our knowledge, for the general instance of K. Therefore, in this work, we develop and formulate coupled oscillator inspired computational models for solving other archetypal combinatorial optimization problems, namely, (a) the Max-K-Cut Problem: defined as the problem of dividing the nodes of a graph G (V, E) into K subsets such that the number of edges across the subsets is maximized (we consider unweighted graphs here); (b) the Graph Coloring Problem: defined as the problem of finding the minimum number of colors required such that every node in the graph is assigned one color, and no two connected nodes (i.e., nodes that share an edge) are assigned the same color; (c) the Traveling Salesman Problem (TSP): which asks the question: "Given a list of cities and the distances between each pair of cities, what is the shortest possible route that visits each city exactly once and returns to the origin city?"; (d) the Hamiltonian Path (and Cycle) Problem: defined as the problem of identifying a path (if it exists) that visits every node in the graph exactly once. The Hamiltonian cycle is the same problem but imposes the additional constraint that the path must return to the node from where it originated; and (e) the Graph Partitioning Problem: defined as the challenge of dividing the vertices of a graph into two equal sets such that the number of common edges is minimized. Building on prior work by Wang *et al.* [10], we develop the Lyapunov functions that map to the problems considered, and subsequently, construct the associated Kuramoto dynamics. We also show that besides phase partitioning/clustering (using external force functions), the relative phase ordering (in appropriately constructed Kuramoto dynamics) can be mapped to the solutions of combinatorial optimization problems such as the traveling salesman problem and the Hamiltonian path problem etc.

Moreover, these example problems also help illustrate the general principle and approach behind developing such models for other problems not covered here. We show that the synchronized oscillator system dynamics can be engineered to solve different optimization problems by appropriately designing the coupling function and the external injection to each oscillator. We also note that in the present work, our focus lies on developing the computational models and not the exact physical implementation. We first discuss the formulation of the Max-K-Cut problem, which can be considered as a general case of the MaxCut problem.

**Max-K-Cut:**

To develop the computational model for the oscillator-inspired Max-K-Cut solver, we first start with the description of the Ising model developed for MaxCut problem and illustrate how the oscillator dynamics can be used to solve the same. The MaxCut problem can be mapped to the Ising Hamiltonian $H$ by representing each node of the corresponding graph by a spin $\{-1, +1\}$. Every edge between the nodes i and j is represented by the interaction coefficient $J_{ij}$ among the spins; $J_{ij} = -1$ when an edge exists between the nodes, and $J_{ij} = 0$ when the edge is absent. Subsequently, solving the MaxCut problem is equivalent to minimizing the Ising Hamiltonian described as:

$$H = - \sum_{i,j=1, i<j}^{N} J_{ij} s_i s_j \tag{1}$$

Wang *et. al* [10] considered the Lyapunov function defined in equation (2) and demonstrated that it is associated with the dynamics of a topologically equivalent (oscillator $\equiv$ spin; coupling element $\equiv$ interaction coefficient; oscillator phase $\equiv$ spin assignment) network of coupled oscillators, under second harmonic injection and is equivalent to $H$.

$$E(\phi(t)) = -C_1 \sum_{i,j=1, j \neq i}^{N} J_{ij} \cos(\Delta\phi_{ij}) - \sum_{i=1}^{N} C_{sync} \cos(2\phi_i(t)) \tag{2}$$

where $C_1$ is the coupling strength, and $C_{sync}$ modulates the strength of the second harmonic injection signal.

Furthermore, equation (2) being a Lyapunov function, the coupled oscillator phases, $\phi_I$, will partition themselves in a way such that the system will tend to minimize $E(\phi(t))$ as the system evolves over time. Consequently, this also minimizes $H$ (Appendix 1). $E(\phi(t))$ is minimized at discrete points $\{0, \pi\}$ in the phase space which effectively correspond to the partitions/sets created by the MaxCut. This enables the oscillator network to directly compute the MaxCut solution as the system evolves towards the ground state energy.

We now extend this approach to the Max-K-Cut problem. To aid the analysis for the Max-K-Cut problem, it is useful to recast the 'spins' that represent each set of the MaxCut as complex quantities, $re^{i\theta_k}$, where $r = 1$, $\theta_k = \frac{2\pi k}{K}$; $k = 1, 2, \ldots, K$. For the MaxCut problem, where $K = 2$, it is evident that the possible spin configurations are $1e^{i\pi} = -1$ ($k = 1$) and $1e^{i(2\pi)} = 1$ ($k = 2$). To map the Max-K-Cut problem, we formulate the "spin" assignment for each of the K sets using the above scheme i.e., each set is assigned a value $re^{i\theta_k}$, where $r = 1$, $\theta_k = \frac{2\pi k}{K}$; $k = 1, 2 \ldots, K$. We note that defining "spins" using this approach for $K > 2$ disengages it from the physical significance of a spin i.e., a spin with an assignment $1e^{i\left(\frac{2\pi}{3}\right)}$ ($k = 1, K = 3$) may not have physical relevance. However, we will continue to use the term here for continuity. It is also noteworthy though that while a complex spin has little physical significance, in oscillator networks, such assignments can be easily represented by the amplitude and phase. In fact, each set created by the Max-K-Cut will be represented by a specific oscillator phase, as discussed further.

The objective function describing the Max-K-Cut (modeled along the same lines as the Ising Hamiltonian) can be described as:

$H_{Max-K-Cut}$

$$= -\sum_{i,j,i<j} J_{ij} \cdot \text{Re}\left(e^{\left(i \lim_{\sigma \to 0} \sum_{k=1}^{K-1}\left(\left((2k-1)\pi - \frac{2k\pi}{K}\right) \cdot e^{-\left(\frac{\left(\Delta\theta_{ij} - \frac{2k\pi}{K}\right)^2}{2\sigma^2}\right)} + \left(\frac{2k\pi}{K} - (2k-1)\pi\right) \cdot e^{-\left(\frac{\left(\Delta\theta_{ij} + \frac{2k\pi}{K}\right)^2}{2\sigma^2}\right)}\right)\right)} s_i s_j^*\right) \quad (3)$$

where $s_i, s_j$ represent the 'spin' assignments given by $re^{i\theta}$; $\Delta\theta_{ij} = \theta_i - \theta_j$, $J_{ij}$ describes the connectivity among the nodes of the graph; $J_{ij} = -1(0)$, when an edge is present (absent) between nodes i and j. For simplicity, we rewrite equation (3) as:

$$H_{Max-K-Cut} = -\sum_{i,j,i<j}^{N} J_{ij} \cdot \text{Re}\left(e^{if(\Delta\theta_{ij})} s_i s_j^*\right) \quad (4)$$

where,

$$f(\Delta\theta_{ij}) = \lim_{\sigma \to 0} \sum_{k=1}^{K-1}\left(\left((2k-1)\pi - \frac{2k\pi}{K}\right) \cdot e^{-\left(\frac{\left(\Delta\theta_{ij} - \frac{2k\pi}{K}\right)^2}{2\sigma^2}\right)} + \left(\frac{2k\pi}{K} - (2k-1)\pi\right) \cdot e^{-\left(\frac{\left(\Delta\theta_{ij} + \frac{2k\pi}{K}\right)^2}{2\sigma^2}\right)}\right) \quad (5)$$

$f(\Delta\theta_{ij})$ is a $2\pi$-periodic odd function. Equation (4) can also be expressed as:

$$H_{\text{Max-K-Cut}} = -\sum_{i,j,i<j}^{N} J_{ij} \cos\left(\Delta\theta_{ij} + f(\Delta\theta_{ij})\right) \tag{6}$$

$H_{\text{Max-K-Cut}}$ is designed such that $-J_{ij}\text{Re}(e^{if(\Delta\theta_{ij})}s_i s_j^*) = -1$ if and only if the nodes corresponding to an edge are in different partitions. For example, consider an edge whose corresponding nodes are placed in the sets with assignments $s_i = 1e^{i\left(\frac{4\pi}{3}\right)}$ (using $k_i = 2$), and $s_j^* = 1e^{-i\left(\frac{2\pi}{3}\right)}$ (using $k_j = 1$), respectively. For this edge, $f(\Delta\theta_{ij}) = \left(\frac{\pi}{3}\right)$; $e^{if(\Delta\theta_{ij})} = e^{i\left(\frac{\pi}{3}\right)}$, and thus, $-J_{ij} \cdot \text{Re}(e^{if(\Delta\theta_{ij})}s_i s_j^*) = -1$. In contrast, if the edge is assigned to the same set (say, using $k_i = k_j = 1$), then $s_i = 1e^{i\left(\frac{2\pi}{3}\right)}$, $s_j^* = 1e^{-i\left(\frac{2\pi}{3}\right)}$ and $f(\Delta\theta_{ij}) = 0$; $e^{i(0)} = 1$. $-J_{ij}\text{Re}(e^{if(\Delta\theta_{ij})}s_i s_j^*)$ then evaluates to 1. Thus, computing the Max-K-Cut of the graph is equivalent to minimizing $H_{\text{Max-K-Cut}}$. Fig. 1 shows the spin assignment and $f(\Delta\theta_{ij})$ for various partitions K.

To emulate the minimization of $H_{\text{Max-K-Cut}}$ in a physical system, we consider a coupled oscillator system with N oscillators. For an oscillator in a coupled network, the phase evolution of that oscillator in the network can be described using the Gen-Adler's equation [25] as described by Wang *et al.* [10],

$$\frac{d\phi_i(t)}{dt} = \omega_i - \omega_{\text{sync}} + \omega_i \sum_{j=1,\ j\neq i}^{N} c_{ij}(\Delta\phi_{ij}) \tag{7}$$

where $\omega_I$ is the natural frequency of the individual oscillator, $\omega_{\text{sync}}$ is the frequency of the synchronized network and $c_{ij}(.)$ is a $2\pi$-periodic function for the coupling among oscillators i and j. The oscillator network is designed to be topologically equivalent to the

| # of Partitions (K) | Spin Assignment | $f(\Delta\theta_{ij})$ |
|---|---|---|
| 2 (Max-Cut) | $1e^{i\pi}(\equiv -1)$<br>$1e^{i(2\pi)}(\equiv 1)$ | |
| 3 | $1e^{i\left(\frac{2\pi}{3}\right)}(\equiv -0.5 + i.\frac{\sqrt{3}}{2})$<br>$1e^{i\left(\frac{4\pi}{3}\right)}(\equiv -0.5 - i.\frac{\sqrt{3}}{2})$<br>$1e^{i\left(\frac{6\pi}{3}\right)}(\equiv 1)$ | |
| 4 | $1e^{i\left(\frac{2\pi}{4}\right)}(\equiv +i)$<br>$1e^{i\left(\frac{4\pi}{4}\right)}(\equiv -1)$<br>$1e^{i\left(\frac{6\pi}{4}\right)}(\equiv -i)$<br>$1e^{i\left(\frac{8\pi}{4}\right)}(\equiv +1)$ | |
| 5 | $1e^{i\left(\frac{2\pi}{5}\right)}(\equiv 0.3 + 0.95i)$<br>$1e^{i\left(\frac{4\pi}{5}\right)}(\equiv -0.8 + 0.58i)$<br>$1e^{i\left(\frac{6\pi}{5}\right)}(\equiv -0.8 - 0.58i)$<br>$1e^{i\left(\frac{8\pi}{5}\right)}(\equiv 0.3 - 0.95i)$<br>$1e^{i\left(\frac{10\pi}{5}\right)}(\equiv 1)$ | |

**Fig. 1.** Spin assignments and corresponding f($\Delta\theta_{ij}$) for various values of K in the Max-K-Cut input graph i.e., oscillator ≡ node; coupling element ≡ edge. Thus, the coupling network can be described by the matrix ($J_{ij}$) of the graph.

Further, assuming that all the oscillators have the same frequency, equal to the synchronized frequency of the network [10], equation (7) evolves to:

$$\frac{d\phi_i(t)}{dt} = -C_1 \sum_{j=1,\ j \neq i}^{N} J_{ij} c_{ij}(\Delta\phi_{ij}) \tag{8}$$

where $C_1$ is a positive constant that signifies the coupling strength among the oscillators.

Under the influence of the injection of the $K^{th}$ harmonic signal to every oscillator in the system, equation (8) further evolves to:

$$\frac{d\phi_i(t)}{dt} = -C_1 \sum_{j=1,\ j \neq i}^{N} J_{ij} c_{ij}(\Delta\phi_{ij}) - C_{sync}\sin(K\phi_i(t)) \tag{9}$$

where $C_{sync}$ is a positive constant that describes the amplitude of the injected signal ($K^{th}$ harmonic). In order to engineer the properties of the oscillator network for solving the Max-K-Cut problem, we carefully design the coupling function $c_{ij}(.)$ to be $\sin\left(\Delta\phi_{ij} + f(\Delta\phi_{ij})\right)$, where $f(\Delta\phi_{ij})$ is derived from equation (5). The system dynamics can then be expressed as:

$$\frac{d\phi_i(t)}{dt} = -C_1 \sum_{j=1,\ j \neq i}^{N} J_{ij} \sin\left(\Delta\phi_{ij} + f(\Delta\phi_{ij})\right) - C_{sync}\sin(K\phi_i(t)) \tag{10}$$

Here, $f(\Delta\phi_{ij})$ can be considered as a coupling function as illustrated in Appendix 2. For a system with the dynamics as described in equation (10), we consider the following Lyapunov function candidate:

$$E(\phi(t)) = -\frac{KC_1}{2} \sum_{i,j,\, j \neq i}^{N} J_{ij} \cos\left(\Delta\phi_{ij} + f(\Delta\phi_{ij})\right) - \sum_{i=1}^{N} C_{sync} \cos(K\phi_i(t)) \qquad (11)$$

To show that $E(\phi(t))$ is a Lyapunov function, we first analyze $\frac{dE(\phi(t))}{dt}$ which can be expressed as $\frac{dE(\phi(t))}{dt} = \frac{\partial E(\phi(t))}{\partial \phi_i(t)} \cdot \frac{d\phi_i(t)}{dt}$. It can be shown (detailed in Appendix 3) that

$$\frac{\partial E(\phi(t))}{\partial \phi_i(t)} = -K \frac{d\phi_i(t)}{dt} \qquad (12)$$

Thus,

$$\frac{dE(\phi(t))}{dt} = \sum_{i=1}^{N} \left[\left(\frac{\partial E(\phi(t))}{\partial \phi_i(t)}\right) \cdot \left(\frac{d\phi_i(t)}{dt}\right)\right] \qquad (13)$$

$$= -K \sum_{i=1}^{N} \left[\left(\frac{d\phi_i(t)}{dt}\right)^2\right] \leq 0 \qquad (14)$$

This indicates that the oscillator network described by the dynamics in equation (10) will aim to minimize the energy function in equation (11) as it evolves over time. $E(\phi(t))$ is minimized at the discrete phase points $\phi_k = \frac{2\pi k}{K}, k = 1, 2, \ldots, K$ which effectively correspond to the partitions / sets ($s_i = re^{i\theta_k}$, where $r = 1$, $\theta_k = \frac{2\pi k}{K}$; $k = 1, 2, \ldots, K$) created by the Max-K-Cut. We note that the minimization of $E(\phi(t))$ at the discrete phase points is facilitated by appropriately designing the injection/forcing function in the dynamics described in equation (9). The role of the forcing function is illustrated in Appendix 4.

Equation (11) is equivalent to $H_{Max-K-Cut}$ in equation (6) with a constant offset which implies that the system will aim to compute the solution to the Max-K-Cut problem as well.

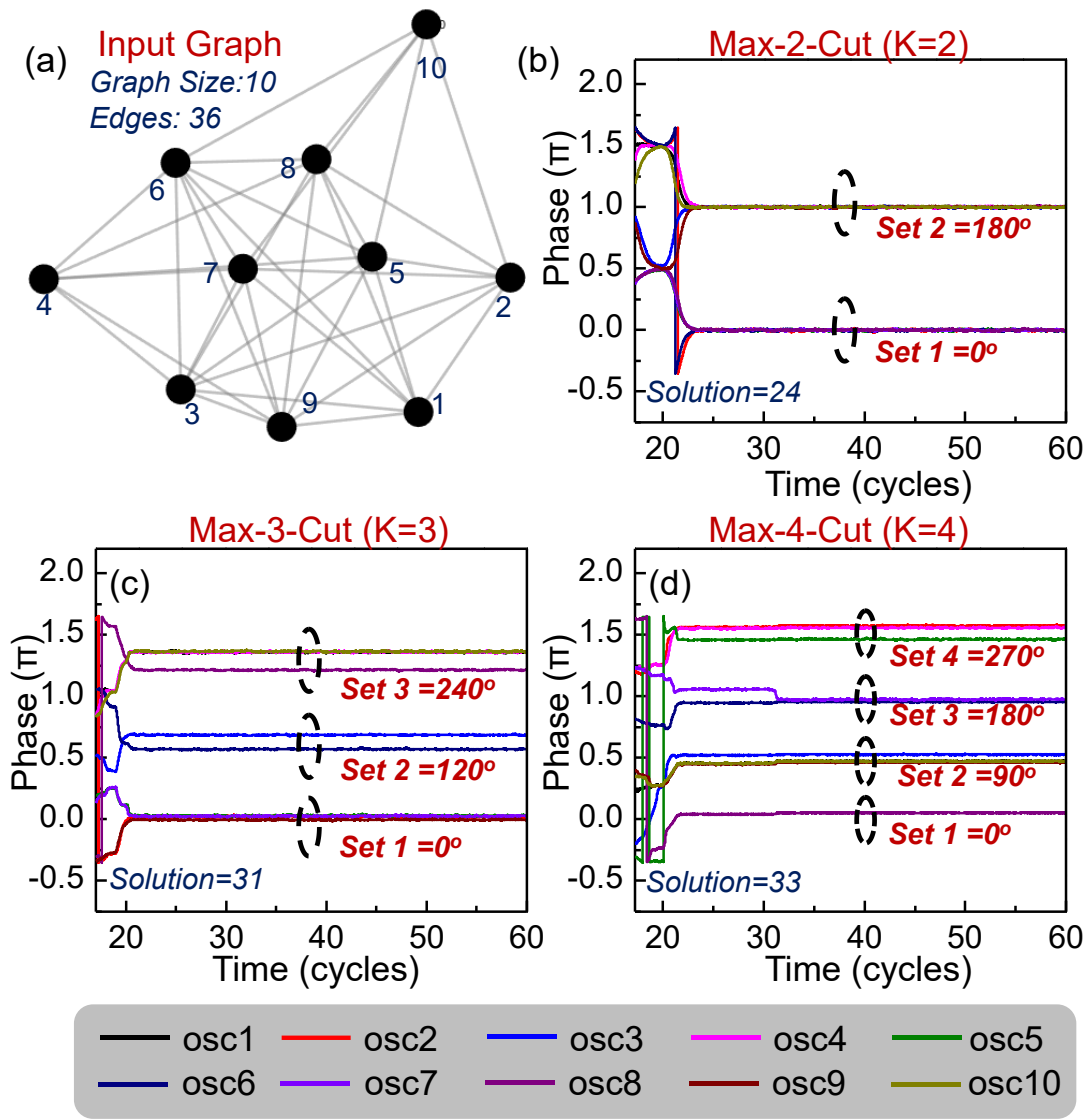

**Fig. 2.** (a) Illustrative 10 node graph; (b)-(d) Evolution of oscillator phases for different values of K (=2,3,4) in the Max-K-Cut problem.

Fig. 2 shows the phase partitions and the resulting Max-K-Cut (K= 2, 3, 4) solutions for an illustrative graph (with 10 nodes) using the oscillator-inspired model developed above. An annealing scheme similar to that used by Wang *et al.* [10], is used to help the system escape from local minima (Appendix 5). It can be observed that the oscillator phases, representing the nodes of the graph, separate into K partitions to yield Max-K-Cut solutions. Further, using the oscillator-inspired solver developed above, we also evaluate some larger graph instances (Appendix 6) from the G-Set dataset for K=2, 3, and 4 partitions. In each case, the system is simulated for 100 cycles and yields high-quality Max-K-Cut solutions.

We next use the above approach to formulate oscillator-inspired computational models for other combinatorial optimization problems.

**Graph Coloring**

Given a graph G(V, E) such that each node is to be assigned a specific color, the graph coloring problem entails finding the minimum number of colors required such that the no two connected nodes (i.e., nodes that share an edge) are assigned the same color. Nodes having the same color can be considered as a set that does not contain edges. Thus, the graph coloring problem can be formulated as finding the smallest K for which the sets created by the Max-K-Cut of the graph do not contain any edges within the set i.e., every edge in the graph are shared between (any) two sets. Solving the problem using this approach entails computing the Max-K-Cut of the graph for different values of K in ascending order, and determining the (smallest) K for which the sets created by the Max-K-Cut does not contain an edge internal to a set. The smallest value of K is the Chromatic number of the graph. This can be considered equivalent to solving the decision version of the graph coloring problem which evaluates if a graph is K-colorable i.e., can a graph be colored with K colors? We note that Crnkic et al. [26] also proposed a similar oscillator-inspired model to solve the graph coloring problem but used a different formulation for the energy function. Some additional constraints on the value of K can be imposed by using graph properties e.g., the value of K will be at most $\Delta+1$, where $\Delta$ is the maximum degree of the graph (Brooks theorem [27]). Fig. 3 shows the coloring for an illustrative example for a graph (same as that considered in Fig. 2) considering K=5 where it can be observed that nodes split into 5 partitions. Moreover, none of the sets contain an edge within the set, implying that the graph is 5 colorable. In contrast, the partitions generated by the Max-4-Cut for the graph (Fig. 2d), which verifies if the graph is 4 colorable, shows that 2 edges lie within the set. Consequently, this implies that the graph is not 4 colorable, and thus, 5 is the minimum number of colors required to color the graph (i.e., chromatic number).

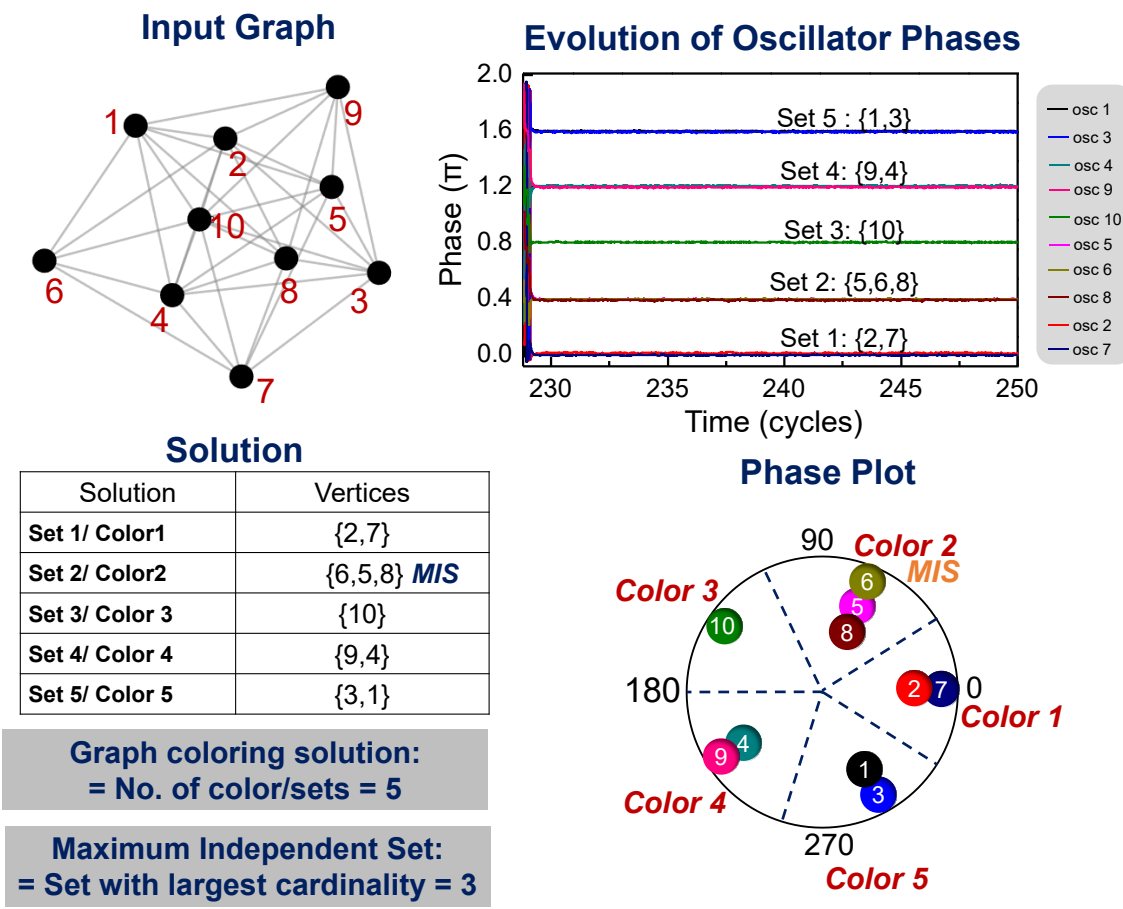

**Fig. 3.** A representative 10 node graph along with the corresponding phase partitions obtained using the oscillator-based computational model (K=5). The graph is 5 colorable since every edge is shared between (any) two sets i.e., there are no edges that lie entirely within one set. Further, the largest such set is the MIS (=3, optimal).

## Maximum Independent Set (MIS) and Maximum Clique (MC)

Each set obtained while computing the graph coloring problem (using the Max-K-Cut formulation) that does not contain an edge corresponds to an independent set of the graph. Thus, the largest independent set among them provides a good approximate solution for the MIS (Fig.3) [28]. Additionally, using the relationship that the Maximum Clique of a graph is the MIS of its complement graph, the same approach can be used to compute the Maximum Clique of a graph by computing the MIS of the complement graph.

## Traveling Salesman Problem (TSP)

Given a list of cities and the distances between each pair, the archetypal TSP is defined as the challenge of finding the shortest possible route that visits each city exactly once and returns to the city of origin. The objective function for the TSP can be formulated as minimizing $H_{TSP}$, where $H_{TSP}$ is defined as:

$$H_{TSP} = - \sum_{i,j,i<j}^{N} J_{ij} \cdot \text{Re}\left(e^{if_{TSP}(\Delta\theta_{ij})} s_i s_j^*\right) \tag{15}$$

where $s_i$ denotes the graph node (city), expressed as a complex quantity described earlier, and $J_{ij}=-D_{ij}$, where $D_{ij}$ refers to the distance between nodes (cities) i and j. Here,

$$f_{TSP}(\Delta\theta_{ij}) = \lim_{\sigma \to 0} - \sum_{\gamma=1, N-1} \left( \frac{2\gamma\pi}{N} \cdot e^{-\left(\frac{\left(\Delta\theta_{ij} - \frac{2\gamma\pi}{N}\right)^2}{2\sigma^2}\right)} + \left(-\frac{2\gamma\pi}{N}\right) \cdot e^{-\left(\frac{\left(\Delta\theta_{ij} + \frac{2\gamma\pi}{N}\right)^2}{2\sigma^2}\right)} \right)$$

$$+ \lim_{\sigma \to 0} \sum_{k=2, k \neq N-1}^{N} \left( \left(\pi - \frac{2k\pi}{N}\right) e^{-\left(\frac{\left(\Delta\theta_{ij} - \frac{2k\pi}{N}\right)^2}{2\sigma^2}\right)} + \left(\frac{2k\pi}{N} - \pi\right) \cdot e^{-\left(\frac{\left(\Delta\theta_{ij} + \frac{2k\pi}{N}\right)^2}{2\sigma^2}\right)} \right) \tag{16}$$

$$H_{TSP} = -\sum_{i,j,i<j}^{N} J_{ij}\cos\left(\Delta\theta_{ij} + f_{TSP}(\Delta\theta_{ij})\right) \qquad (17)$$

The TSP is equivalent to computing a phase ordering such that the sum of distances ($D_{ij}$) among the adjacent nodes (cities) in the phase ordering is as small as possible. To achieve this, $f_{TSP}(\Delta\theta_{ij})$ is designed such that $-J_{ij}\text{Re}(e^{if_{TSP}(\Delta\theta_{ij})}s_i s_j^*) = D_{ij}$ if the nodes are adjacent (i.e., $\text{Re}(e^{if_{TSP}(\Delta\theta_{ij})}s_i s_j^*)=1$), else $-J_{ij}\text{Re}(e^{if_{TSP}(\Delta\theta_{ij})}s_i s_j^*) = -D_{ij}$ (i.e., $\text{Re}(e^{if_{TSP}(\Delta\theta_{ij})}s_i s_j^*) = -1$); in other words, the system incurs a larger energy penalty for putting nodes with larger distances adjacent to each other. Therefore, minimizing $H_{TSP}$ minimizes the sum of distances among the adjacent nodes. Furthermore, the circular phase ordering ensures that dynamics not only compute the shortest route among the cities but also guarantee that one returns back to the original city after visiting each city only once. Along similar lines as the Max-K-Cut problem, the corresponding oscillator dynamics can be expressed as:

$$\frac{d\phi_i(t)}{dt} = -C_1 \sum_{j=1,j\neq i}^{N} J_{ij}\sin\left(\Delta\phi_{ij} + f_{TSP}(\Delta\phi_{ij})\right) - C_{sync}\sin(N\phi_i(t)) \qquad (18)$$

As discussed earlier, the force term on the right-hand side ensures that the oscillators settle to a phase value of $\frac{2\pi k}{K}; k = 1, 2, \ldots, N$. Similar to the Max-K-Cut analysis, the Lyapunov function candidate considered for these dynamics is:

$$E(\phi(t)) = -\frac{NC_1}{2} \sum_{i,j,\, j\neq i}^{N} J_{ij}\cos\left(\Delta\phi_{ij} + f_{TSP}(\Delta\phi_{ij})\right) - \sum_{i=1}^{N} C_{sync}\cos(N\phi_i(t)) \qquad (19)$$

Which can be shown to be minimized over time. Further, equation (19) is equivalent to $H_{TSP}$ (with a constant offset), and thus, will be minimized over time. Fig. 4 shows a representative example for a 10 city TSP graph solved using the coupled oscillators.

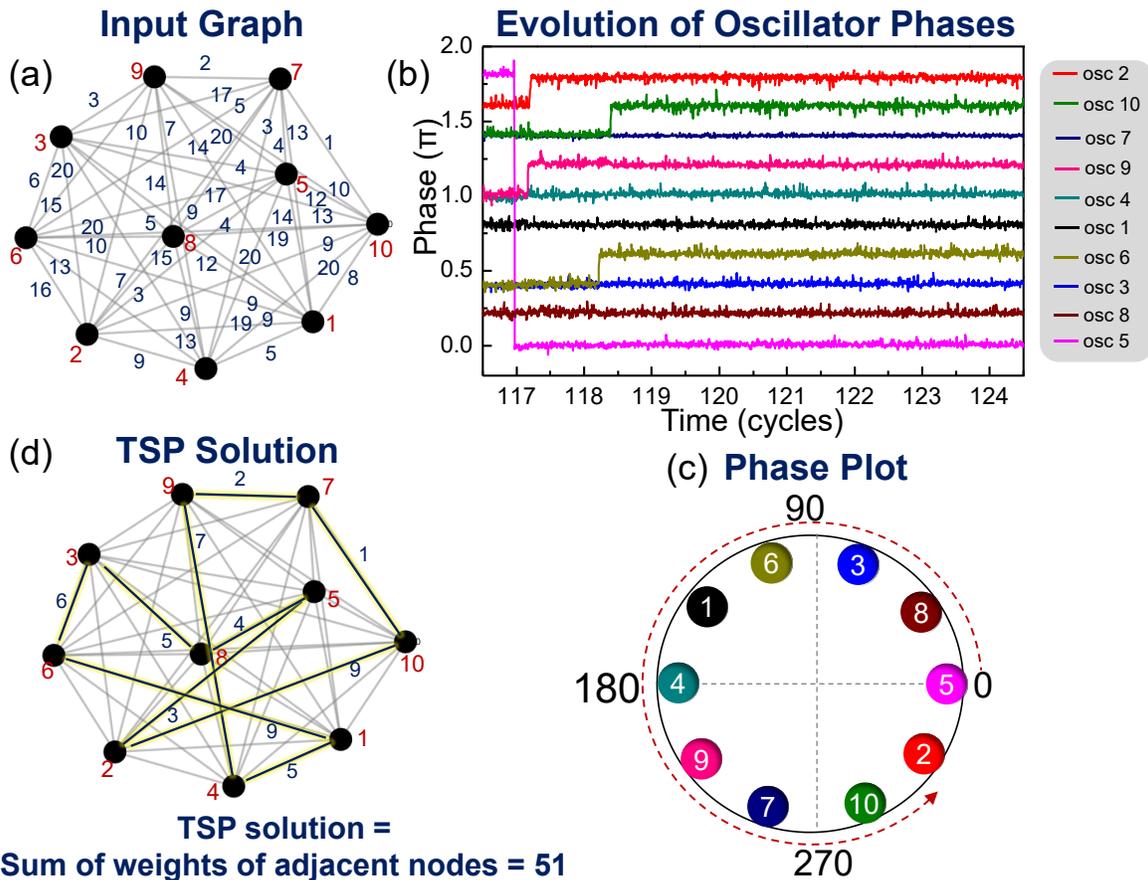

**Fig. 4.** (a) A representative 10 node graph; (b) oscillator phase evolution in a topologically equivalent oscillator system. (c) phase plot and the corresponding phase ordering that encodes the TSP solution. (d) TSP solution obtained from the phase ordering (=51; optimal for the above graph).

**Hamiltonian Cycle/Hamiltonian Path**

The Hamiltonian path problem is defined as the problem of identifying a path (if it exists) that visits every node exactly once; computing the Hamiltonian cycle entails solving the same problem with the added constraint of returning back to the node of origin. In terms of the oscillator system dynamics, the objective here is to design a topologically equivalent coupled oscillator network that yields a phase ordering such that the number of edges among the adjacent nodes is maximized. Subsequently, if every node shares an edge with the adjacent nodes in the phase ordering, the resulting ordering represents the Hamiltonian cycle (and path); if only one edge is missing, then only a Hamiltonian path (but not cycle) exists. The Hamiltonian for the above problem can be formulated using the same approach used for TSP with the exception that the system is 'rewarded' (in terms of energy) for bringing connected nodes with edges adjacent to each other. The objective function for the Hamiltonian cycle path can be expressed as:

$$H_{HC} = -\sum_{i,j,i<j}^{N} J_{ij} \cdot \text{Re}\left(e^{if_{HC}(\Delta\theta_{ij})} s_i s_j^*\right) \tag{20}$$

where,

$$f_{HC}(\Delta\theta_{ij}) = \lim_{\sigma \to 0} -\sum_{\gamma=1,N-1} \left( \left(\pi - \frac{2\gamma\pi}{N}\right) \cdot e^{-\left(\frac{\left(\Delta\theta_{ij} - \frac{2\gamma\pi}{N}\right)^2}{2\sigma^2}\right)} + \left(\frac{2\gamma\pi}{N} - \pi\right) \cdot e^{-\left(\frac{\left(\Delta\theta_{ij} + \frac{2\gamma\pi}{N}\right)^2}{2\sigma^2}\right)} \right)$$

$$+ \lim_{\sigma \to 0} \sum_{k=2,k \neq N-1}^{N} \left( \left(\frac{\pi}{2} - \frac{2\pi k}{N}\right) e^{-\left(\frac{\left(\Delta\theta_{ij} - \frac{2k\pi}{N}\right)^2}{2\sigma^2}\right)} + \left(\frac{2\pi k}{N} - \frac{\pi}{2}\right) \cdot e^{-\left(\frac{\left(\Delta\theta_{ij} + \frac{2k\pi}{N}\right)^2}{2\sigma^2}\right)} \right) \tag{21}$$

$$H_{HC} = -\sum_{i,j,i<j}^{N} J_{ij} \cos\left(\Delta\theta_{ij} + f_{HC}(\Delta\theta_{ij})\right) \tag{22}$$

It can be observed that $f_{HC}(\Delta\theta_{ij})$ has a similar form to $f_{TSP}(\Delta\theta_{ij})$ defined for TSP (equation (16)) with the exception that $\Delta\theta_{ij} + f_{HC}(\Delta\theta_{ij})$ converges to $\pi$ when a pair of connected nodes are adjacent to each other (i.e., $\Delta\theta_{ij} = \pm\frac{2\pi}{N}$) in the phase space whereas $\Delta\theta_{ij} + f_{TSP}(\Delta\theta_{ij})$ converges to π when they are in non-adjacent positions (i.e., $\Delta\theta_{ij} = \pm\frac{2k\pi}{N}; 1 < k \leq N; k \neq 1, N-1$). In these respective (desired) scenarios, the system is 'rewarded' by lowering energy (equations (17) and (22)). Further, $\Delta\theta_{ij} + f_{HC}(\Delta\theta_{ij})$ (for the Hamiltonian cycle /path problem) converges to $\frac{\pi}{2}$ when the connected nodes are in non-adjacent positions (undesired case) whereas $\Delta\theta_{ij} + f_{TSP}(\Delta\theta_{ij})$ (for the TSP) converges to 0 when the nodes are in adjacent to each other in the phase space. In the case of the Hamiltonian cycle/path problem, this implies that the system does not lower energy when the connected nodes are in non-adjacent positions (i.e., system is rewarded only for placing nodes adjacent to each other) whereas in the case of TSP the system energy actually increases when the connected nodes are placed next to each other in the phase space. We note that similar to the TSP, $\Delta\theta_{ij} + f_{HC}(\Delta\theta_{ij})$ can be designed to converge to 0 (instead of $\frac{\pi}{2}$) by modifying the pre-factors for the gaussian distributions in the second term in the right-hand side of equation (21). This would impose an energy penalty similar to TSP when a pair of connected nodes attains an undesirable configuration. Our motivation in choosing $\frac{\pi}{2}$ here was to illustrate how different energy functions can 'designed'. We note that in the both the cases (i.e., $\Delta\theta_{ij} + f_{HC}(\Delta\theta_{ij})$ converging to $\frac{\pi}{2}$ or 0, when a pair of connected nodes are non-adjacent to each other), the equivalence between the optimal solution to the Hamiltonian cycle / path problem and the ground state of the system does not change, and the system will still continue to evolve towards a lower energy state in both the cases.

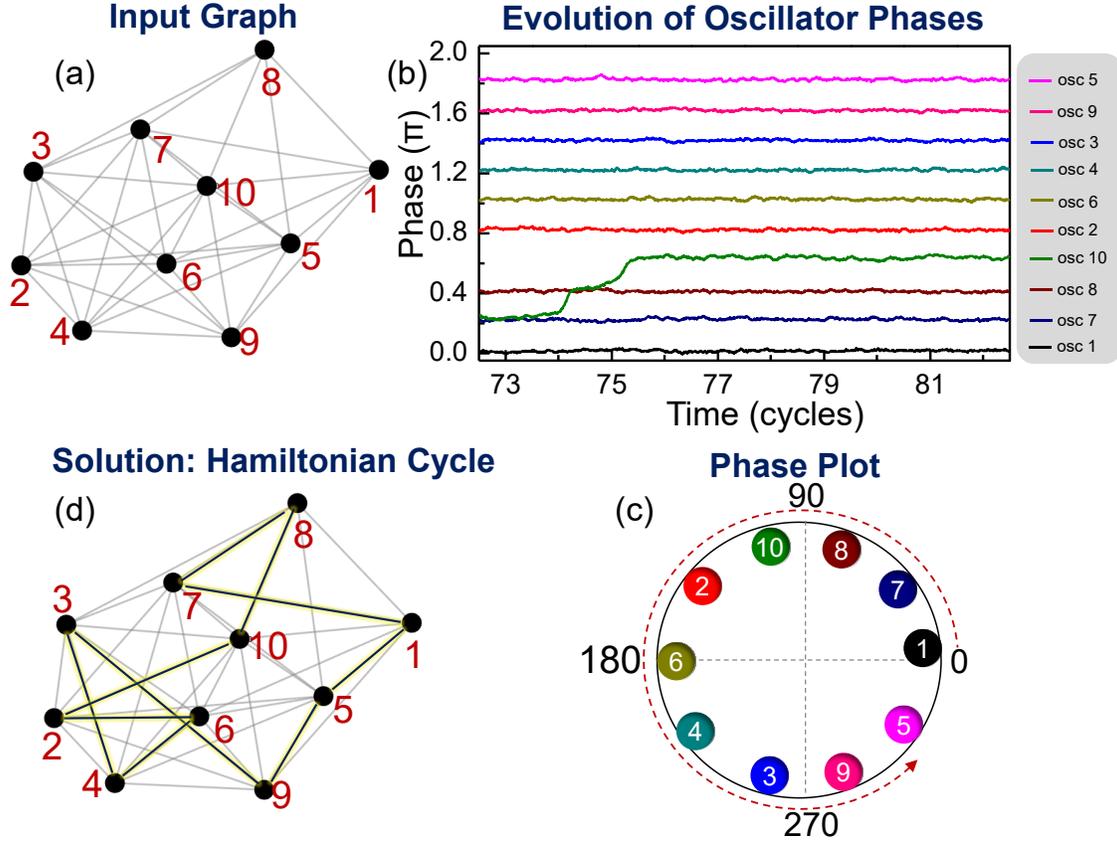

**Fig. 5.** (a) A representative 10 node randomly instantiated graph; (b) phase dynamics of the equivalent coupled oscillator network for solving the Hamiltonian path/cycle problem; (c) corresponding ordering of phases, resulting in a Hamiltonian cycle as shown in (d).

Similar to the TSP, the corresponding oscillator dynamics and the energy function of the system can be described by,

$$\frac{d\phi_i(t)}{dt} = -C_1 \sum_{j=1, j \neq i}^{N} J_{ij} \sin\left(\Delta\phi_{ij} + f_{HC}(\Delta\phi_{ij})\right) - C_{\text{sync}} \sin(N\phi_i(t)) \qquad (23)$$

$$E(\phi(t)) = -\frac{NC_1}{2} \sum_{i,j,\, j \neq i}^{N} J_{ij} \cos\left(\Delta\phi_{ij} + f_{HC}(\Delta\phi_{ij})\right) - \sum_{i=1}^{N} C_{\text{sync}} \cos(N\phi_i(t)) \qquad (24)$$

Equation (24), which is equivalent to $H_{HC}$ in equation (20) (with an offset), is minimized over time and in the process evolves towards the solution of the Hamiltonian cycle. Fig. 5 shows the Hamiltonian cycle computed on a demonstrative problem using the above approach.

**Graph Partitioning**

The problem of dividing the nodes of a graph into two equal sets such that the number of common edges is minimized. For simplicity, we assume here that the number of nodes (N) is even. The objective function for the problem can be defined as:

$$H_{GP} = -\sum_{i,j,i<j}^{N} J_{ij}\left(A.\text{Re}\left(\left(\prod_{n=1}^{N} s_n\right)^{\frac{1}{N}} \cdot e^{-i\frac{\pi}{2}}\right) + B.\text{Re}(s_i s_j^*)\right) \quad (25)$$

Here, $s_j = 1e^{i\theta_j}$, where $\theta_j \in \{0, \pi\}$. The first term on the RHS in equation (20) determines the relative reward (energy reduction) for dividing the set into two equal parts, and the

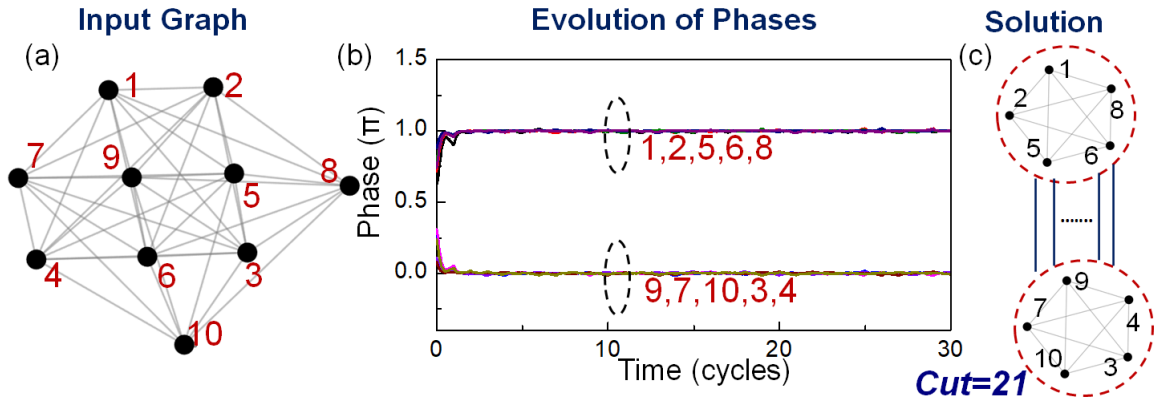

**Fig. 6.** (a) A representative 10 node graph and (b) phase dynamics of the coupled oscillators for the representative graph. (c) Corresponding graph partition solution of 21 (optimal).

second term determines the reward for minimizing the number of common edges. A (>0)

and B (>0) are constants that determine the ratio of the two rewards and determine which constraint has preponderance in the optimization process. Further, $J_{ij} = +1(0)$ when an edge is present (absent).

To solve the problem using coupled oscillators, the objective can be expressed as designing coupled oscillator network where the oscillators exhibit a phase bi-partition $\{0, \pi\}$ such that each set has an equal number of oscillators, and the number of edges among the two clusters is small as possible. These properties can be realized using a system that exhibits the following dynamics:

$$\frac{d\phi_i(t)}{dt} = -\sum_{j=1, j\neq i}^{N} J_{ij} \left( 2C_1 \sin\left(\frac{1}{N}\sum_{n=1}^{N} \phi_n - \frac{\pi}{2}\right) + C_2 \sin(\Delta\phi_{ij}) \right) - C_{sync} \sin(2\phi_i(t)) \quad (26)$$

with the corresponding Lyapunov function:

$$E(\phi(t)) = -\sum_{i,j=1, j\neq i}^{N} J_{ij} \left( 2NC_1 \cos\left(\frac{1}{N}\sum_{n=1}^{N} \phi_n - \frac{\pi}{2}\right) + C_2 \cos(\Delta\phi_{ij}) \right) - \sum_{i=1}^{N} C_{sync} \cos(2\phi_i(t)) \quad (27)$$

The solution to the graph partitioning problem then corresponds to the two sets of oscillators corresponding to a phase value of 0 and π, respectively, as shown in Fig. 6.

**Conclusion**

This work elucidates how the natural energy minimization in oscillator-based dynamical systems can be used to develop computational models for solving combinatorial optimization problems. While the above models develop the relationships between the oscillator dynamics and the objective function characterizing the optimization problem, we note that local energy minima in the energy function, the precision of the weights etc.

will impact the quality of the calculated solution and the time-to-compute, and thus, the overall success of the approach. Additionally, the effect of the strength of oscillator coupling as well as the external injection on the phase clustering, the role of chimera synchronization [29-35], and its impact on the dynamical and computational properties also needs to be investigated further. We would like to point out that classical implementations of the above models are unlikely to reduce the fundamental computational hardness of the above problems although many recent works [21-24,36] have shown that empirically, such analog models and approaches could provide a substantial speed up with potential trade-offs with other metrics such as solution accuracy; a detailed understanding of the scaling behavior of such analog system-inspired computational models is of increasing interest and needs to be advanced further [37-39]. However, the models developed in this work help expand the applications of this computing approach beyond the traditional combinatorial problems (e.g., MaxCut) considered in prior works, and thus, help bolster the case for exploring analog inspired non-Boolean approaches to computationally intractable optimization problems.

**Appendix 1:**

Here, we illustrate that minimizing the Lyapunov function for the MaxCut problem in equation (2) is equivalent to minimizing H in equation (1), first shown by Wang *et al.* [10]. Similar analysis can be extended to the other combinatorial optimization problems considered here.

The Lyapunov function for the oscillator system to solve the MaxCut is given by:

$$E(\phi(t)) = -C_1 \sum_{i,j=1, j\neq i}^{N} J_{ij} \cos(\Delta\phi_{ij}) - \sum_{i=1}^{N} C_{sync} \cos(2\phi_i(t)) \quad (A1.1)$$

where, $C_1 > 0$, $C_{sync} > 0$; $J_{ij} = -1$.

Equation (A1.1) (same as equation (2)) will achieve a minimum when $\cos(\Delta\phi_{ij}) = -1$ (resulting from the first term on the right-hand side) and $\cos(2\phi_i(t)) = 1$ (resulting from the second term on the right-hand side of the equation).

These conditions are satisfied (facilitating E to achieve a minimum) when $\phi_i = 0$ or $\pi$. At these discrete phase points, equation (A1.1) can be expressed as:

$$E(\phi(t)) = -C_1 \sum_{i,j=1, j\neq i}^{N} J_{ij} \cos(\Delta\phi_{ij}) - NC_{sync} \quad (A1.2)$$

since $\cos(2\phi_i(t))=1$ for $\phi_i=0$ or $\pi$.

Further, by mapping every spin $s_i = +1$ to $\phi_i = 0$ ($\pi$), and $s_i = -1$ to $\phi_i = \pi$ (0), it can be observed that $\cos(\Delta\phi_{ij}) = -1$, when the nodes lie in the opposite set and $\cos(\Delta\phi_{ij}) = +1$ when the nodes lie in the same set. This is similar to $s_i s_j = -1$ (+1) when the spins have opposite (same) alignment, respectively. Furthermore, the phase configuration $[\phi_1, \phi_2, \phi_3, \ldots \phi_N]$ (where $\phi_i = 0$ or $\pi$) which minimizes $E(\phi(t))$ will correspond to the spin configuration $[s_1, s_2, s_3 \ldots, s_N]$ that minimizes the Ising Hamiltonian $H$ given by:

$$H = -\sum_{i,j=1, i<j}^{N} J_{ij} s_i s_j \tag{A1.3}$$

In fact, with $C_1 = \frac{1}{2}$, equation (A1.2) is equivalent to the Ising Hamiltonian with a constant offset.

$$E(\phi(t)) = -\sum_{i,j=1,i<j}^{N} J_{ij} \cos(\Delta\phi_{ij}) - NC_{sync} \tag{A1.4}$$

Here, $\frac{1}{2}$ is used since the spin interaction in the Hamiltonian is represented by two terms ($J_{ij}$, $J_{ji}$) in the Lyapunov function.

## Appendix 2:

Here, we illustrate how $f(\Delta\phi_{ij})$ can essentially be considered as a coupling function. Equation (10) can be expanded and rewritten as:

$$\frac{d\phi_i(t)}{dt} = -C_1 \sum_{j=1, i \neq j}^{N} \left[ \mathbf{J_{ij}cos}\left(\mathbf{f(\Delta\phi_{ij})}\right) \cdot \sin(\Delta\phi_{ij}) + \left(\mathbf{J_{ij}sin}\left(\mathbf{f(\Delta\phi_{ij})}\right)\right) \cdot \cos(\Delta\phi_{ij}) \right] \\ - C_{sync}\sin(K\phi_i) \tag{A2.1}$$

Here, $J_{ij}\cos\left(f(\Delta\phi_{ij})\right)$ and $J_{ij}\sin\left(f(\Delta\phi_{ij})\right)$ can be considered as coupling functions between the quadrature outputs of the oscillator.

**Appendix 3:** Proof to show $\frac{\partial E(\phi(t))}{\partial \phi_i(t)} = -K \frac{d\phi_i(t)}{dt}$

$$E(\phi(t)) = -\frac{KC_1}{2} \sum_{i,j,\,j \neq i}^{N} J_{ij} \cos\left(\Delta\phi_{ij} + f(\Delta\phi_{ij})\right) - \sum_{i=1}^{N} C_{sync} \cos(K\phi_i(t)) \quad (A3.1)$$

$$\frac{\partial E(\phi(t))}{\partial \phi_i(t)} = -\frac{KC_1}{2} \sum_{j=1,\,j \neq i}^{N} J_{ij} \frac{\partial}{\partial \phi_i(t)} \cos\left(\phi_i - \phi_j + f(\phi_i - \phi_j)\right)$$

$$-\frac{KC_1}{2} \sum_{j=1,\,j \neq i}^{N} J_{ji} \frac{\partial}{\partial \phi_i(t)} \cos\left(\phi_j - \phi_i + f(\phi_j - \phi_i)\right) - C_{sync} \frac{\partial}{\partial \phi_i(t)} \cos(K\phi_i) \quad (A3.2)$$

$$\frac{\partial E(\phi(t))}{\partial \phi_i(t)} = C_2 \sum_{j=1,\,j \neq i}^{N} J_{ij} \sin\left(\Delta\phi_{ij} + f(\Delta\phi_{ij})\right) \cdot \left(1 + \frac{\partial}{\partial \phi_i(t)} f(\Delta\phi_{ij})\right)$$

$$+ C_2 \sum_{j=1,\,j \neq i}^{N} J_{ji} \sin\left(\Delta\phi_{ji} + f(\Delta\phi_{ji})\right) \cdot \left(-1 + \frac{\partial}{\partial \phi_i(t)} f(\Delta\phi_{ji})\right) + C_{sync} K \cdot \sin(K\phi_i) \quad (A3.3)$$

$$\frac{\partial}{\partial \phi_i(t)} f(\Delta\phi_{ij}) = \lim_{\sigma \to 0} \sum_{k=1}^{K-1} \left( -\left((2k-1)\pi - \frac{2k\pi}{K}\right) \cdot e^{-\left(\frac{\left(\phi_i - \phi_j - \frac{2k\pi}{K}\right)^2}{2\sigma^2}\right)} \cdot \left(\frac{2\left(\phi_i - \phi_j - \frac{2k\pi}{K}\right)}{2\sigma^2}\right) \right.$$

$$\left. - \left(\frac{2k\pi}{K} - (2k-1)\pi\right) \cdot e^{-\left(\frac{\left(\phi_i - \phi_j + \frac{2k\pi}{K}\right)^2}{2\sigma^2}\right)} \cdot \left(\frac{2\left(\phi_i - \phi_j + \frac{2k\pi}{K}\right)}{2\sigma^2}\right) \right) \quad (A3.4)$$

$$\frac{\partial}{\partial \phi_i(t)} f(\Delta\phi_{ji}) = \lim_{\sigma \to 0} \sum_{k=1}^{K-1} \left((2k-1)\pi - \frac{2k\pi}{K}\right) \cdot e^{-\left(\frac{\left(\phi_j - \phi_i - \frac{2k\pi}{K}\right)^2}{2\sigma^2}\right)} \cdot \left(\frac{2\left(\phi_j - \phi_i - \frac{2k\pi}{K}\right)}{2\sigma^2}\right) \quad \text{(A3.5)}$$

$$+ \left(\frac{2k\pi}{K} - (2k-1)\pi\right) \cdot e^{-\left(\frac{\left(\phi_j - \phi_i + \frac{2k\pi}{K}\right)^2}{2\sigma^2}\right)} \cdot \left(\frac{2\left(\phi_j - \phi_i + \frac{2k\pi}{K}\right)}{2\sigma^2}\right)$$

Using the relation that, $\lim_{\sigma \to 0} \frac{e^{-\alpha^2/\sigma^2}}{\sigma^2} = 0$ in equations (A3.4) and (A3.5)

$$\frac{\partial}{\partial \phi_i(t)} f(\Delta\phi_{ij}) = \frac{\partial}{\partial \phi_i(t)} f(\Delta\phi_{ji}) = 0 \quad \text{(A3.6)}$$

Similarly, $\frac{\partial}{\partial \phi_i(t)} f_{TSP}(\Delta\phi_{ij}) = \frac{\partial}{\partial \phi_i(t)} f_{TSP}(\Delta\phi_{ji})$ relevant to solving the TSP can be shown to be equal to 0 as illustrated further on.

Substituting equation (A3.6) into (A3.3) and using $\sin(x) = -\sin(-x)$ & $J_{ij} = J_{ji}$ we get,

$$\frac{\partial E(\phi(t))}{\partial \phi_i(t)} = 2C_2 \sum_{j=1, j\neq i}^{N} J_{ij} \sin\left(\Delta\phi_{ij} + f(\Delta\phi_{ij})\right) + C_{sync} K.\sin(K\phi_i) \quad \text{(A3.7)}$$

$$\frac{\partial E(\phi(t))}{\partial \phi_i(t)} = 2 \cdot \frac{K.C_1}{2} \sum_{j=1, j\neq i}^{N} J_{ij} \sin\left(\Delta\phi_{ij} + f(\Delta\phi_{ij})\right) + C_{sync} K.\sin(K\phi_i) \quad \text{(A3.8)}$$

$$\frac{\partial E(\phi(t))}{\partial \phi_i(t)} = KC_1 \sum_{j=1, j\neq i}^{N} J_{ij} \sin\left(\Delta\phi_{ij} + f(\Delta\phi_{ij})\right) + C_{sync} K.\sin(K\phi_i) = -K.\frac{d\phi_i(t)}{dt} \quad \text{(A3.9)}$$

$$\frac{dE(\phi(t))}{dt} = \sum_{i=1}^{N}\left[\left(\frac{\partial E(\phi(t))}{\partial \phi_i(t)}\right)\cdot\left(\frac{d\phi_i(t)}{dt}\right)\right] = -K\sum_{i=1}^{N}\left[\left(\frac{d\phi_i(t)}{dt}\right)^2\right] \leq 0 \qquad (A3.10)$$

To show $\frac{\partial}{\partial \phi_i(t)} f_{TSP}(\Delta\phi_{ij}) = \frac{\partial}{\partial \phi_i(t)} f_{TSP}(\Delta\phi_{ji}) = 0$.

$$\frac{\partial}{\partial \phi_i(t)} f_{TSP}(\Delta\phi_{ij})$$

$$= \lim_{\sigma \to 0} - \sum_{\gamma=1,N-1} \left( \left(-\frac{2\gamma\pi}{N}\right)\cdot e^{-\left(\frac{\left(\phi_i - \phi_j - \frac{2\gamma\pi}{N}\right)^2}{2\sigma^2}\right)} \cdot \left(\frac{2\left(\phi_i - \phi_j - \frac{2\gamma\pi}{N}\right)}{2\sigma^2}\right) \right.$$

$$\left. + \left(\frac{2\gamma\pi}{N}\right)\cdot e^{-\left(\frac{\left(\phi_i - \phi_j + \frac{2\gamma\pi}{N}\right)^2}{2\sigma^2}\right)} \cdot \left(\frac{2\left(\phi_i - \phi_j + \frac{2\gamma\pi}{N}\right)}{2\sigma^2}\right) \right) \qquad (A3.11)$$

$$+ \lim_{\sigma \to 0} \sum_{k=2;k\neq N-1}^{N} \left( -\left(\pi - \frac{2k\pi}{N}\right)\cdot e^{-\left(\frac{\left(\phi_i - \phi_j - \frac{2k\pi}{N}\right)^2}{2\sigma^2}\right)} \cdot \left(\frac{2\left(\phi_i - \phi_j - \frac{2k\pi}{N}\right)}{2\sigma^2}\right) \right.$$

$$\left. - \left(\frac{2k\pi}{N} - \pi\right)\cdot e^{-\left(\frac{\left(\phi_i - \phi_j + \frac{2k\pi}{N}\right)^2}{2\sigma^2}\right)} \cdot \left(\frac{2\left(\phi_i - \phi_j + \frac{2k\pi}{N}\right)}{2\sigma^2}\right) \right)$$

$$\frac{\partial}{\partial \phi_i(t)} f_{TSP}(\Delta\phi_{ji}) \quad \text{(A3.12)}$$

$$= \lim_{\sigma \to 0} -\sum_{\gamma=1,N-1} \left( \left(\frac{2\gamma\pi}{N}\right) \cdot e^{-\left(\frac{\left(\phi_j - \phi_i - \frac{2\gamma\pi}{N}\right)^2}{2\sigma^2}\right)} \cdot \left(\frac{2\left(\phi_j - \phi_i - \frac{2\gamma\pi}{N}\right)}{2\sigma^2}\right) \right.$$

$$\left. + \left(-\frac{2\gamma\pi}{N}\right) \cdot e^{-\left(\frac{\left(\phi_j - \phi_i + \frac{2\gamma\pi}{N}\right)^2}{2\sigma^2}\right)} \cdot \left(\frac{2\left(\phi_j - \phi_i + \frac{2\gamma\pi}{N}\right)}{2\sigma^2}\right) \right)$$

$$+ \lim_{\sigma \to 0} \sum_{k=2}^{N} \left( \left(\pi - \frac{2k\pi}{N}\right) \cdot e^{-\left(\frac{\left(\phi_j - \phi_i - \frac{2k\pi}{N}\right)^2}{2\sigma^2}\right)} \cdot \left(\frac{2\left(\phi_j - \phi_i - \frac{2k\pi}{N}\right)}{2\sigma^2}\right) \right.$$

$$\left. + \left(\frac{2k\pi}{N} - \pi\right) \cdot e^{-\left(\frac{\left(\phi_j - \phi_i + \frac{2k\pi}{N}\right)^2}{2\sigma^2}\right)} \cdot \left(\frac{2\left(\phi_j - \phi_i + \frac{2k\pi}{N}\right)}{2\sigma^2}\right) \right)$$

Using the relation that, $\lim_{\sigma \to 0} \frac{e^{-\alpha^2/\sigma^2}}{\sigma^2} = 0$ in equations (A3.11) and (A3.12)

$$\frac{\partial}{\partial \phi_i(t)} f_{TSP}(\Delta\phi_{ij}) = \frac{\partial}{\partial \phi_i(t)} f_{TSP}(\Delta\phi_{ji}) = 0 \quad \text{(A3.13)}$$

Substituting equation (A3.13) into (A3.3), and using a similar approach as described above, it can be shown that the energy function for the TSP (equation (19) in main text) is a Lyapunov function. The approach can be used to show that $\frac{\partial}{\partial \phi_i(t)} f_{HC}(\Delta\phi_{ij}) = \frac{\partial}{\partial \phi_i(t)} f_{HC}(\Delta\phi_{ji}) = 0$, and the corresponding energy function for the Hamiltonian cycle/path (equation (24) in the main text) is a Lyapunov function.

**Appendix 4:**

Here, we elucidate the role of the injection term/ external forcing function in discretizing the phases for the general case of the Max-K-Cut problem. We first consider the system dynamics for solving the Max-K-Cut problem *without the external injection*:

$$\frac{d\phi_i(t)}{dt} = -C_1 \sum_{j=1,\ j\neq i}^{N} J_{ij} \sin\left(\Delta\phi_{ij} + f(\Delta\phi_{ij})\right) \qquad (A4.1)$$

The corresponding Lyapunov function for this system is given by:

$$E(\phi(t)) = -\frac{KC_1}{2} \sum_{i,j,\ j\neq i}^{N} J_{ij} \cos\left(\Delta\phi_{ij} + f(\Delta\phi_{ij})\right) \qquad (A4.2)$$

Equation (A4.2) achieves a minimum when the phase difference is $\frac{2\pi k}{K}$, $k = 1, 2, \ldots, K$. We also note that $f(\Delta\phi_{ij})$ in the first term is specifically designed such that the $\cos\left(\Delta\phi_{ij} + f(\Delta\phi_{ij})\right)$ term equals $-1$ (minimum value for the cos(.) function) when $\Delta\phi_{ij} = \frac{2\pi k}{K}$ (for MaxCut, $f(\Delta\phi_{ij})=0$ since the cos(.) naturally achieves a minimum at a phase difference of π).

While this minimum is attained when $\Delta\phi_{ij} = \frac{2\pi k}{K}$, $\phi_i$ and $\phi_j$ can assume analog values. For instance, in the case of the Max-3-Cut problem, $\phi_i = 0$ and $\phi_j = \frac{2\pi}{3}$, and $\phi_i = \frac{\pi}{3}$ and $\phi_j = \pi$ are both probable solutions since $\Delta\phi_{ij} = \frac{2\pi}{3}$ in both cases ($f(\Delta\phi_{ij}) = \frac{\pi}{3}$ in both cases). Consequently, the oscillators exhibit a continuous distribution of phases, that cannot be directly mapped to the K-sets created by the Max-K-Cut.

In contrast, when the force term is considered, the resulting system dynamics are described by,

$$\frac{d\phi_i(t)}{dt} = -C_1 \sum_{j=1,\ j\neq i}^{N} J_{ij} \sin\left(\Delta\phi_{ij} + f(\Delta\phi_{ij})\right) - C_{sync}\sin(K\phi_i(t)) \quad \text{(A4.3)}$$

for which the corresponding Lyapunov function is given by,

$$E(\phi(t)) = -\frac{KC_1}{2} \sum_{i,j,\ j\neq i}^{N} J_{ij} \cos\left(\Delta\phi_{ij} + f(\Delta\phi_{ij})\right) - \sum_{i=1}^{N} C_{sync}\cos(K\phi_i(t)) \quad \text{(A4.4)}$$

In order for Equation (A4.4) to achieve a minimum, not only does $\Delta\phi_{ij} + f(\Delta\phi_{ij}) = (2m + 1)\pi$ has to hold true (condition arises from the first term), but $\phi_i$ and $\phi_j$ have to also equal $\frac{2\pi k}{K}$. In other words, the forcing term discretizes the phases to K points. Subsequently, these phase points can be directly mapped to the K sets created by the Max-K-Cut problem.

**Appendix 5:**

The strength of the coupling among the oscillators ($C_1$) is increased linearly given by $C_1 = 1 + t\frac{A-1}{T}$; where T is the total simulation time (equivalent to 100 oscillation cycles), and A is a constant which is chosen such that the maximum value of $C_1 = A$ (at $t = T$). This implementation is inspired by prior work by Wang et al. [10], wherein increasing $C_1$ linearly was effectively used as an annealing schedule that helped the system better escape from local minima, and thus, facilitated improved solution quality. We observe similar behavior as shown in the representative graph (Fig. 7), wherein increasing the coupling strength linearly (while maintaining $C_1$ below a critical threshold) produced better solutions compared to maintaining a constant coupling strength.

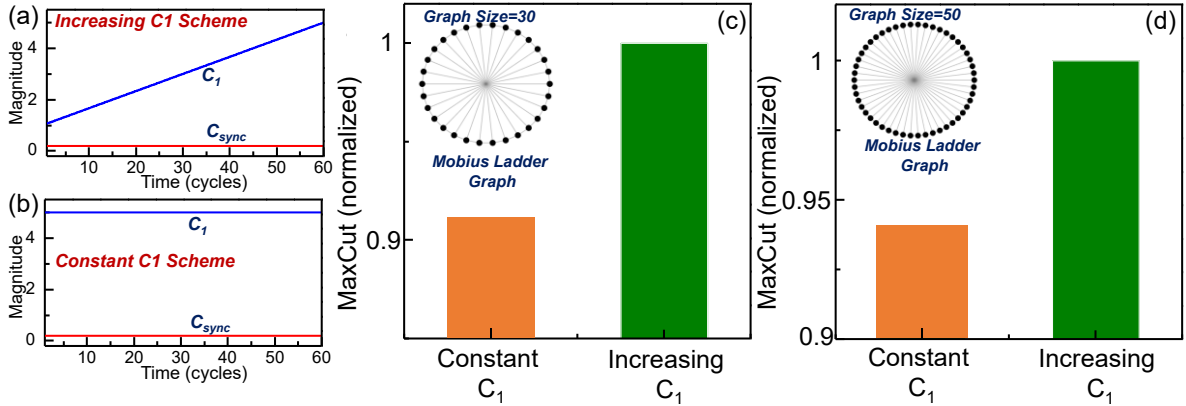

**Fig. 7.** Impact of $C_1$ evolution on the MaxCut solution quality. (a)(b) Two schemes for evolution of the coupling constants considered in the example. $C_1$ increases linearly in (a), while it is constant in (b); $C_{sync}$ is constant in both the cases. (c)(d) Observed MaxCut solutions (normalized to optimal solution) for two representative Mobius ladder graphs of size 30 and 50 nodes, respectively. It can be observed empirically that linearly increasing $C_1$ produces improved solutions.

However, it must be noted that considering the negative nature of coupling, $C_1$ cannot be made indefinitely large and must be below a certain threshold. In fact, we observe that the $\frac{C_1}{C_{sync}}$ ratio (i.e., coupling strength among the oscillator nodes relative to the coupling constant for external forcing function) is critical to achieving the desired functionality, and the system clusters get destabilized if the $\frac{C_1}{C_{sync}}$ ratio is above a certain threshold. This is illustrated for a representative graph in Fig. 8 where it can be observed that the clusters start to destabilize for $\frac{C_1}{C_{sync}} > 70$. We note that the exact value of the threshold will depend on the properties of the graph such as size and connectivity as well as the computing problem.

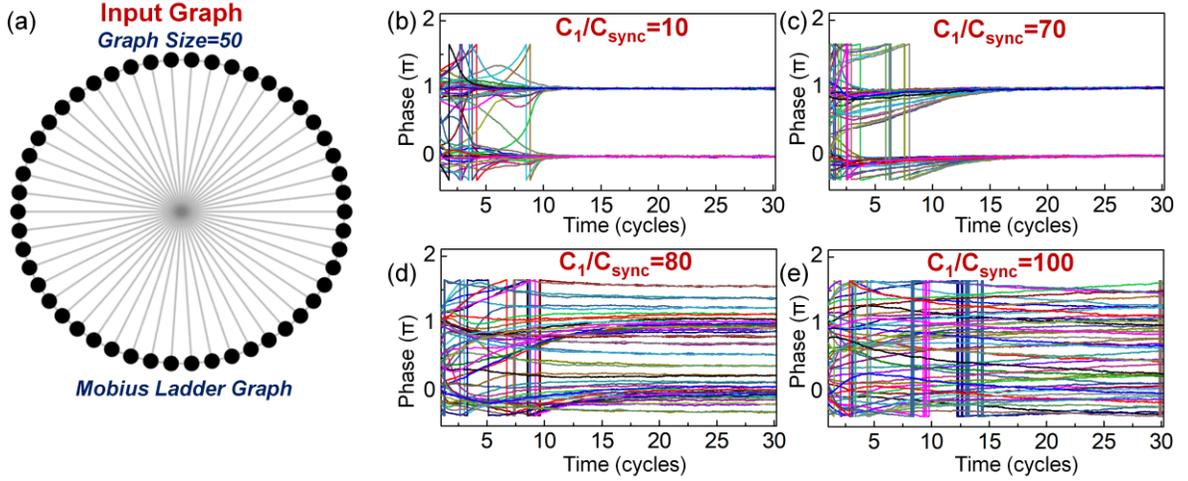

**Fig. 8.** (a) Mobius ladder graph considered in the illustrative example. (b)-(e) Evolution of oscillator phases for different $C_1/C_{sync}$ ratios (=10, 70, 80, 100, respectively). $C_1$: coupling strength among the oscillators; $C_{sync}$: Coupling strength of external injection signal / force function. It can be observed that $C_1$ (relative to $C_{sync}$) must remain below a critical threshold in order to observe the phase clustering.

We also note that besides $C_1, C_{sync}$ the standard deviation σ of the gaussian distributions must be carefully designed for the problems considered in this work. While force function ensures that certain angles $\left(\frac{2k\pi}{N}; k = 1, 2, 3..N\right)$ are energetically favored, the system, while evolving towards this low energy configuration, may assume other phase angles not equal to $\frac{2k\pi}{N}$. It is important to make sure that such angles do not help the system evolve towards a lower energy configuration. This can be ensured by appropriately choosing σ as well as adding additional gaussian terms in the function $f(\Delta\theta_{ij})$ designed to penalize the system (i.e., increase energy) for such configurations. For instance, in the case of the Hamiltonian cycle, gaussian functions can be added to ensure that the system energy is increased for all angles not equal to the $\pm\frac{2\pi}{N}$ (instead of just considering an energy penalty for $\frac{2k\pi}{N}; k \neq 1, N-1$.

**Appendix 6:**

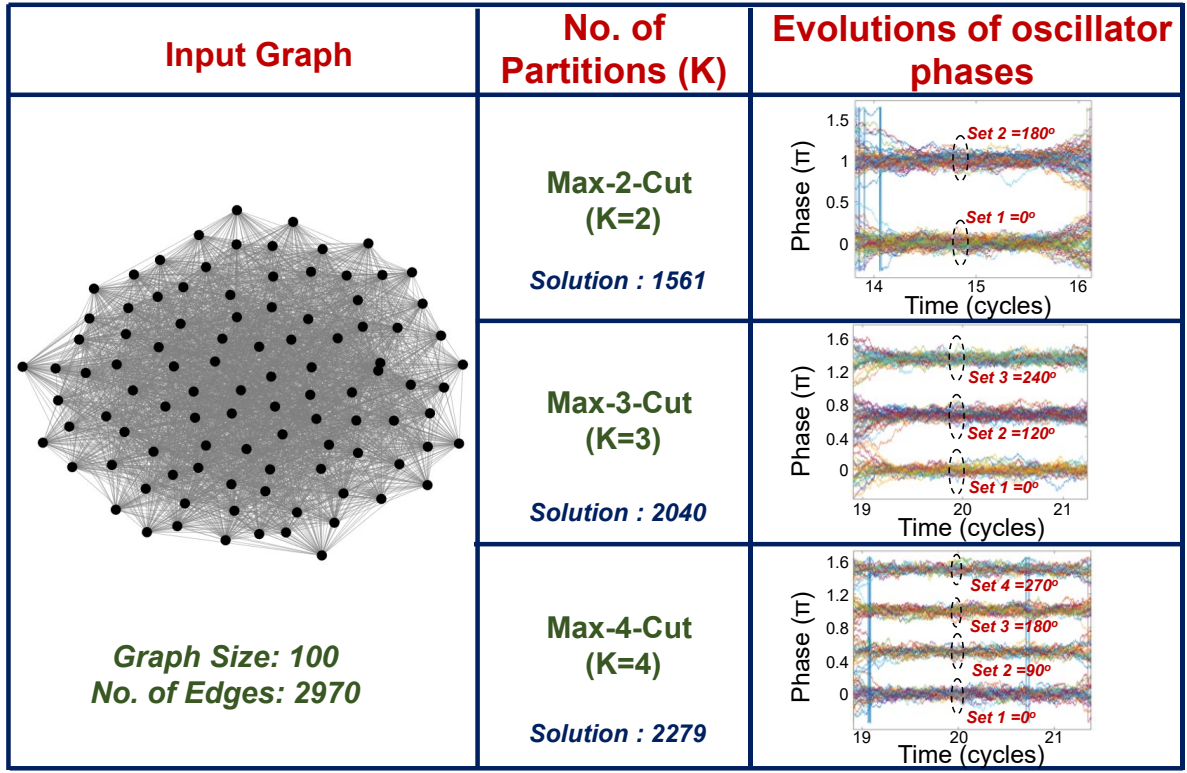

**Fig. 9.** A representative 100 node (randomly instantiated) graph and the evolution of oscillator phases corresponding to the solution of the Max-K-Cut problem (K=2, 3, 4).

| Graph | Nodes (V) | Edges (E) | Oscillator Solution (Max-2-cut) | Solution Quality (Max-2-cut) | Oscillator Solution (Max-3-cut) | Solution Quality (Max-3-cut) | Oscillator Solution (Max-4-cut) | Solution Quality (Max-4-cut) |
|---|---|---|---|---|---|---|---|---|
| G1 | 800 | 19176 | 11624 | 100% | 15032 | 99.1% | 16166 | 96.2% |
| G2 | 800 | 19176 | 11606 | 99.8% | 14878 | 98.1% | 16197 | 96.4% |
| G3 | 800 | 19176 | 11613 | 99.9% | 14852 | 97.88% | 16204 | 96.4% |
| G4 | 800 | 19176 | 11635 | 99.9% | 14894 | 98.1% | 16208 | 96.4% |
| G5 | 800 | 19176 | 11570 | 99.4% | 14915 | 98.2% | 16211 | 96.4% |

**Fig. 10.** Graph instances from the G-Set solved using the oscillator-based computational model; the oscillators compute yields high quality Max-*K*-Cut solutions within >96% of the solutions obtained by F. Ma et al. [40]. Here, Solution Quality $= \frac{\text{OScillator Solution}}{\text{Best Known Solution}} \times 100\%$.

**Acknowledgment:**

This work was supported in part by NSF ASCENT grant (No. 2132918).